\documentclass[12pt,final,a4paper,reqno,psamsfonts]{amsart}
\usepackage[latin1]{inputenc}
\usepackage{amssymb,amscd,amsxtra,verbatim,ifthen,graphicx}
\usepackage[english]{babel}
\usepackage[all]{xy}
\usepackage{hyperref}
\usepackage{enumitem}

\mathchardef\ordinarycolon\mathcode`\: \mathcode`\:=\string"8000
\begingroup \catcode`\:=\active
  \gdef:{\mathrel{\mathop\ordinarycolon}}
\endgroup

\theoremstyle{plain}
\newtheorem{theorem}{Theorem}
\newtheorem{proposition}[theorem]{Proposition}

\newtheorem{lemma}[theorem]{Lemma}

\newtheorem{corollary}[theorem]{Corollary}

\theoremstyle{definition}
\newtheorem{definition}[theorem]{Definition}

\theoremstyle{remark}
\newtheorem{remark}[theorem]{Remark}

\newcommand{\pr}{\mathrm{pr}}

\newcommand{\spn}{\mathrm{span}}

\newcommand{\supp}{\mathrm{supp}}

\newcommand{\reg}{\mathrm{reg}}

\newcommand{\Ad}{\mathrm{Ad}}

\newcommand{\copol}{\mathrm{copol}}
\newcommand{\cohom}{\mathrm{cohom}}

\newcommand{\diag}{\mathrm{diag}}

\newcommand{\vol}{\mathrm{vol}}
\newcommand{\Pu}{\mathrm{Pu}}
\newcommand{\zero}{\mathrm{Zero}}

\newcommand{\lie}{\mathrm{Lie}}
\newcommand{\tr}{\mathrm{tr}}

\renewcommand{\Re}{\mathrm{Re}}

\newcommand{\On}{\mathbf{O}}
\newcommand{\on}{\mathfrak{o}}
\newcommand{\Son}{\mathbf{SO}}

\newcommand{\Un}{\mathbf{U}}
\newcommand{\un}{\mathfrak{u}}
\newcommand{\Sun}{\mathbf{SU}}

\newcommand{\Spn}{\mathbf{Sp}}
\renewcommand{\spn}{\mathfrak{sp}}

\newcommand{\iso}{\mathrm{Iso}}

\newcommand{\R}{\mathbf{R}}
\newcommand{\C}{\mathbf{C}}
\renewcommand{\H}{\mathbf{H}}

\newcommand{\N}{\mathbf{N}}
\newcommand{\K}{\mathbf{K}}

\newdir^{ (}{{}*!/-5pt/@^{(}}
\newcommand{\ka}{{\mathcal A}}

\newcommand{\kc}{{\mathcal C}}
\newcommand{\kd}{{\mathcal D}}
\newcommand{\ke}{{\mathcal E}}

\newcommand{\kh}{{\mathcal H}}

\newcommand{\ks}{{\mathcal S}}

\newcommand{\gothg}{{\mathfrak g}}
\newcommand{\gothh}{{\mathfrak h}}

\newcommand{\gothm}{{\mathfrak m}}
\newcommand{\gothn}{{\mathfrak n}}

\newcommand{\gothS}{{\mathfrak S}}

\begin{document}
\title[An Integration Formula for Isometric Actions]{A general Weyl-type Integration Formula\\ for Isometric Group Actions}
\author{Frederick Magata}
\maketitle


\begin{abstract}
We show that integration over a $G$-manifold $M$ can be reduced to integration over a \emph{minimal section} $\Sigma$ with respect to an induced weighted measure and integration over a homogeneous space $G\!/\!N$. We relate our formula to integration formul\ae\ for \emph{polar actions} and calculate some weight functions. In case of a compact Lie group acting on itself via conjugation, we obtain a classical result of Hermann Weyl. Our formula allows to view almost arbitrary isometric group actions as \emph{generalized random matrix ensembles}. We also establish a reductive decomposition of Killing fields with respect to a minimal section.
\end{abstract}

\section{Introduction and main results}
A particularly ``nice'' class of isometric group actions are \emph{polar} actions.
In this situation there is a complete, connected and embedded
submanifold $\Sigma\subseteq M$, which intersects with every orbit
and which is perpendicular to them in the intersection points. Such
a $\Sigma$ is then called a \emph{section} for the polar action.
Examples arise from the actions of compact Lie groups on themselves
via conjugation or the isotropy actions of Riemannian symmetric
spaces. In these cases sections are given by maximal tori, resp.
maximal flat submanifolds. Polar actions have been extensively studied in the literature (see e.g. \cite{BCO, PT2} for results and references).
What makes them nice is the fact that their orbit geometry is very special. For instance, principal orbits of polar representations are isoparametric submanifolds. Also, the orbit space structure is considerably simple: it is isometric to an orbifold as a path metric space. Many of the properties one observes for maximal tori in compact Lie groups (e.g. any two maximal tori are conjugate, Chevalley's restriction theorem, etc.) also hold for sections of polar actions.

In \cite{GOT} Gorodski, Olmos and Tojeiro introduced a
generalization of the concept of a section in order to measure the
defect of an arbitrary isometric action from being polar. A \textbf{fat section} $\Sigma$ (or \emph{$k$-section} in loc. cit.) is defined
for an isometric action of a Lie group $G$ on a Riemannian
manifold $M$ as follows:
\begin{enumerate}[label=(\Alph{*}), ref=(\Alph{*})]
\item\label{pr:1} $\Sigma$ is a complete, connected, embedded and totally
geodesic submanifold of $M$,

\item\label{pr:2} $\Sigma$ intersects every orbit of the $G$-action,

\item\label{pr:3} for all $G$-regular $p\in \Sigma$ we have $\nu_p(G\cdot
p)\subseteq T_p\Sigma$ with codimension $k$,

\item\label{pr:4} for all $G$-regular $p\in \Sigma$ and $g\in G$ such that $g\cdot
p\in \Sigma$ we have $g\cdot \Sigma=\Sigma$.
\end{enumerate}
The sections of a polar action are $0$-sections in the above sense. And vice versa, if an isometric
action allows $0$-sections, then it is in fact polar. This motivates
us to call the minimal integer $k$, for which a $k$-section exists, the \textbf{copolarity}
$\copol(G,M)$ of the isometric action $(G,M)$. By a \textbf{minimal section} we then
mean any $\copol(G,M)$-section. One checks easily, that $M$ is always a $k$-section with $k$ being the dimension of a principal orbit in this case. Also, certain components of the fixed point set of a principal isotropy group give examples of fat sections (see \cite{GOT} or \cite{Mag1} for details).
Particularly interesting is the case when the minimal sections are different from $M$. We then say that the action has \textbf{non-trivial} copolarity. For an effective action, this is the case whenever the principal isotropy groups are non-trivial.

Although the definition appears to be technical, fat sections have a nice structure theory, which in fact is similar to
the one of sections for polar actions. By structure theory we
mean that we associate to each fat section $\Sigma$ the Lie
group $W:=N_G(\Sigma)\!/\!Z_G(\Sigma)$ (the normalizer of $\Sigma$ in
$G$ modulo its pointwise stabilizer), which acts on $\Sigma$ in a
natural way. This $W$ is called the \textbf{fat Weyl group} of
$\Sigma$, although in general it can be almost any Lie group and
hence is not a Weyl group in the classical sense. There are strong
relations between the action of $G$ on $M$ and the one of $W$ on
$\Sigma$. Basically, all information concerning the orbit space structure of the group action $(G,M)$ can be obtained from the group action $(W,\Sigma)$. This follows from the fact that $G\backslash M$ and $W\backslash\Sigma$ are canonically isometric as path metric spaces. This immediately gives a $\kc^0$-version of Chevalley's restriction theorem\begin{footnote}{We do not know whether the $\kc^\infty$-version of this theorem is true in general, or not.}\end{footnote}. Also, any two minimal sections are conjugate to each other (i.e. they can be mapped by the group action onto each other). Note that if $\Sigma$ is a minimal section, then $W$ is $\copol(G,M)$-dimensional. See \cite{GOT} or \cite{Mag1} for further results on fat sections and their structure theory.

In this paper we also deal with the relation between the action of $G$
on $M$ and the one of $W$ on $\Sigma$. Our first result is that if we fix a minimal section $\Sigma$, then the $G$-Killing fields can be decomposed reductively into those which coincide with $W$-Killing fields on $\Sigma$ along $\Sigma$ and those which are everywhere perpendicular to $\Sigma$ along $\Sigma$.
\begin{theorem}\label{t:reductive_decomp}
Let $\Sigma$ be a minimal section of $(G,M)$ and $N=N_G(\Sigma)$, $H=Z_G(\Sigma)$.
We identify the $G$-Killing fields on $M$ with $\gothg=\lie(G)$. Let further $\gothn=\lie(N),\ \gothh=\lie(H)$ and $\gothm=\{X\in\gothg\mid X|_\Sigma \text{ is everywhere perpendicular to } \Sigma\}$. Then
$$\gothg=\gothn+\gothm \text{ and } \gothh=\gothn\cap\gothm.$$

This decomposition is reductive in the sense that $[\gothn,\gothm]\!\subseteq\!\gothm$. Even better, 
\mbox{$\Ad_g(\gothm)\subseteq\gothm$} holds for all $g\in N$. Thus $\gothg/\gothh=\gothn/\gothh\oplus\gothm/\gothh$ is an $\Ad_G(N)$-invariant direct sum.
\end{theorem}
\begin{corollary}
If $G$ above is connected, then $G$ is generated by $N^\circ$ and $\exp(\gothm)$. The subgroup $K\le G$ generated by $\exp(\gothm)$ is normal in $G$ and $G$ is a quotient of $N^\circ\ltimes K$.
\end{corollary}

In the situation of Theorem \ref{t:reductive_decomp}, we call a left-$G$- and right-$W$-invariant Riemannian metric on $G\!/\!H$ \textbf{$(G$-$W)$-invariant} ($W$ acts naturally in a proper and free fashion from the right on $G\!/\!H$). A $(G$-$W)$-invariant metric is called \textbf{adapted}, if the decomposition $\gothg/\gothh=\gothn/\gothh\oplus\gothm/\gothh$ is orthogonal. See Section \ref{s:decomp} for further remarks. Our main result is:
\begin{theorem}[Weyl's integration formula]\label{t:integration}
Under the assumptions of Theorem \ref{t:reductive_decomp} let $W=W(\Sigma)=N\!/\!H$. We further assume that
$G\!/\!H$ carries an adapted $(G$-$W)$-invariant metric and that $W\subseteq G\!/\!H$ carries the induced metric and $G\!/\!N$ the quotient metric via $W\hookrightarrow G\!/\!H\twoheadrightarrow G\!/\!N$. Furthermore all manifolds are equipped with their corresponding Riemannian measure. Then:
\begin{enumerate}
\item For every $f\in L^1(M)$
$$\int_M f(x)\, dx=\int_{G\!/\!N} \left(\int_\Sigma f(g\cdot s)\, \delta_\ke(s)ds\right)\, d(gN).$$

\item For every $f\in L^1(M)^G$
$$\int_M f(x)\, dx=\vol(G\!/\!N)\int_\Sigma f(s)\,\delta_\ke(s)ds.$$
In particular, if $L^1(M)^G\ne\{0\}$, then $\vol(G\!/\!N)$ is finite.

\item If $G\!/\!N$ is compact, the assignment
$$\Theta_p:L^p(M)^G\to L^p(\Sigma)^W,\, f\mapsto\sqrt[p]{\vol(G\!/\!N)\delta_\ke}\,f|_\Sigma$$
is a surjective linear isometry for any $1\le p<\infty$.

\item If $G\!/\!H$ has finite volume, then so do $G\cdot s, W\cdot s, W$ and $G\!/\!N$. Furthermore, $\delta_\ke(s)$ is, up to a covering factor, the ratio of two volume scaling factors:
$$\delta_\ke(s)=\left\{\begin{array}{cl}
0 & \text{if } s \text{ is singular,}\\
\frac{|G_s\!/\!H|}{|W_s|}\cdot \frac{{\vol(G\cdot s)}/{\vol(G\!/\!H)}}{{\vol(W\cdot s)}/{\vol(W)}}& \text{if } s \text{ is regular or exceptional}.
\end{array}\right.
$$
\end{enumerate}
\end{theorem}
Although (iv) in the theorem above already reveals the geometrical nature of the function $\delta_\ke$, its actual definition is slightly more technical. Let $\omega_s:G\!/\!H\to G\cdot s$ denote the ``orbit map'' of $s\in\Sigma$. This $\omega_s$ is a diffeomorphism, a finite covering or a fiber bundle, depending on whether $G\cdot s$ is a principal, exceptional or singular orbit. For this we remark that since $\Sigma$ is a minimal section and $H=Z_G(\Sigma)$ it follows that $H$ is the isotropy group of a principal orbit and hence all isotropy groups of principal orbits coincide along $\Sigma$. The differential of $\omega_s$ in $eH$ restricted to $\gothm/\gothh$ is the map
$$d\omega_s(eH)|_{\gothm/\gothh}:\gothm/\gothh\to\nu_s\Sigma,\ X+\gothh\mapsto X_s.$$
Now $\delta_\ke$ is defined on $\Sigma$ by
$$\delta_\ke(s):=|\det(d\omega_s(eH)|_{\gothm/\gothh})|.$$
\begin{proposition}\label{p:delta} The function $\delta_\ke$ is non-negative and contained in $\kc(\Sigma)^W\cap\kc^\infty(\Sigma^\reg)$. It vanishes exactly in the $G$-singular points. It has a continuous, resp. smooth, $G$-invariant continuation to $M$, resp. $M^\reg$, also denote by $\delta_\ke$.
\end{proposition}

If the action $(G,M)$ is polar, then $G\!/\!H\twoheadrightarrow G\!/\!N$ is a $|W|$-fold covering and we obtain
\begin{corollary}\label{c:polar} Under the assumptions of Theorem \ref{t:integration} let $(G,M)$ be polar. Then
\begin{enumerate}
\item If $G\!/\!H$ has finite volume, then $W$ is finite. Furthermore, $\delta_\ke(s)$ is, up to a covering factor, a volume scaling function:
$$\delta_\ke(s)=\left\{\begin{array}{cl}
0 & \text{if } s \text{ is singular,}\\
|G_s\!/\!H|\cdot \frac{\vol(G\cdot s)}{\vol(G\!/\!H)}& \text{if } s \text{ is regular or exceptional}.
\end{array}\right.
$$
\item If $G\!/\!H$ has finite volume, then for every $f\in L^1(M)$
\begin{eqnarray*}\int_M f(x)\, dx&=&\frac{1}{|W|}\int_{G\!/\!H} \left(\int_\Sigma f(g\cdot s)\, \delta_\ke(s)ds\right)\, d(gH)\\
&=&\frac{1}{|W|}\int_\Sigma \left(\int_{G\!/\!H}f(g\cdot s)d(gH)\right)\, \delta_\ke(s)ds.
\end{eqnarray*}
\item If $G\!/\!H$ has finite volume, then for every $f\in L^1(M)^G$
\begin{eqnarray*}\int_M f(x)\, dx&=&\tfrac{\vol(G\!/\!H)}{|W|}\int_\Sigma f(s)\,\delta_\ke(s)ds\\
&=&\tfrac{1}{|W|}\int_\Sigma f(s)\,\vol(G\cdot s)ds.
\end{eqnarray*}
\end{enumerate}
\end{corollary}
In case that $G$ is a compact Lie group with bi-invariant metric, acting on itself via conjugation, then any maximal torus $T$ is both $\Sigma$ and $H$. Hence, formula (iii) above is the classical integration formula of Weyl in this case. Corollary \ref{c:polar} has been independently discovered by several authors (e.g. \cite{AWY1,AWY2,F,GT,Mag}). The authors of the first three papers do not mention the notion of a polar action, however. 

In \cite{AWY1} one can find the definition of a \emph{quasi-smooth measure} used in

\begin{corollary}
Under the assumptions of Theorem \ref{t:integration}, if $\vol(M)<\infty$, then $G\!/\!N$ and $\Sigma$ have finite volume, too. The latter with respect to the quasi-smooth measure $\delta_\ke\, ds$:
$$\vol(G\!/\!N)=\tfrac{\vol(M)}{\vol_\ke(\Sigma)}, \text{where } \vol_\ke(\Sigma)=\int_\Sigma\delta_\ke(s)\, ds.$$
\end{corollary}

\begin{remark}
With the notation and assumptions as in Theorem \ref{t:integration}. Let $p\in\kc^\infty(M)^G$. Inspired by \cite{AWY1} we call $(G,\varphi,M,p(x)dx,W(\Sigma),\Sigma,ds)$ a \emph{generalized random matrix ensemble}. Furthermore we call $M$ the \emph{integration manifold}, $\Sigma$ the \emph{eigenvalue manifold with symmetry $W$}, $p(s)\delta_\ke(s)\,ds$ the \emph{generalized eigenvalue distribution} and $p(s)\delta_\ke(s)$ the \emph{generalized joint density function}. We stress, that we do not assume the action to be polar here. We hope that this point of view, associating generalized random matrix ensembles with almost arbitrary isometric actions, may be useful in the context of mathematical physics.
\end{remark}

In some cases $\delta_\ke(s)$, resp. $\vol(G\cdot s)$, can be calculated explicitly. For polar actions related with symmetric spaces and so called \emph{Hermann actions}, this can be accomplished by some ``root space decomposition'', (see \cite{DK} or \cite{H2} and \cite{GT} for explicit formulas). Below we compute $\delta_\ke$ for the $k$-fold direct sums of the standard representations of $\Son(n), \Sun(n)$ and $\Spn(n)$, which for $k\ge 2$ are examples of non-polar actions:
\begin{proposition}\label{p:direct_sum}
Let $\K$ denote one of the skew fields $\R, \C$ or $\H$ and denote the standard representation of $G=\Son(n), \Sun(n)$ or $\Spn(n)$ on $\K^n$ by $\rho_n, \mu_n$, resp. $\nu_n$. The copolarity of the $k$-fold direct sum representations of $\rho_n, \mu_n$ or $\nu_n$ is non-trivial if and only if $1\le k\le n-1$. More precisely, besides the polar case $k=1$, we have for $2\le k\le n-1$:
$$\begin{array}{|c|c|c|c|c|c|c|c|}
\hline \varphi      & H         & \cohom(G,V)       & \Sigma     & \copol(G,V)        & W(\Sigma) \\ \hline\hline
       k\cdot\rho_n & \Son(n-k) & \frac{k(k+1)}{2}  & \R^{k^2}   & \frac{k(k-1)}{2}   & \On(k)    \\ \hline
       k\cdot\mu_n  & \Sun(n-k) & k^2               & \R^{2k^2}  & k^2                & \Un(k)    \\ \hline
       k\cdot\nu_n  & \Spn(n-k) & k(2k-1)           & \R^{4k^2}  & k(2k+1)            & \Spn(k)   \\ \hline
\end{array}$$
Here $H$ is a principal isotropy group along $\Sigma$. In each case, $\varphi$ is left-matrix multiplication of elements of $G$ on $V=\K^{n\times k}$ (as a real space). $H$ sits in $G$ as block matrices $\diag(\mathbf{1},A), A\in H$ 
and elements of $\Sigma$ have the form $p=(B,\mathbf{0})^t, B\in\K^{k\times k}$.
$G\!/\!H$ admits a $(G$-$W)$-invariant Riemannian metric. Let $d:=\dim_\R(\K)$ and $p$ be as above. Then
$$\delta_\ke(p)=\frac{1}{\sqrt{2^{dk(n-k)}}}|\det(B)|^{d(n-k)}.$$
For $\K=\H$, the determinant is to be understood in the sense of Dieudonné \cite{Die}. 
\end{proposition}
In \cite{GOT} it is proved that the irreducible taut representations are exactly those whose copolarity is less than or equal to one. Gorodski showed in \cite{Gor} that the representations appearing in Proposition \ref{p:direct_sum} are all \emph{taut}. So Proposition \ref{p:direct_sum} yields further examples for the fact that reducible taut representations cannot be characterized by the before mentioned criterion.

The rest of the paper is organized as follows: In Section 2 we prove the decomposition Theorem and characterize adapted metrics. Section 3 is devoted to the proof of Weyl's integration formula and its corollaries. In Section 4 we prove Proposition \ref{p:direct_sum}.

\section{Decomposition of Killing Fields and Adapted Metrics}\label{s:decomp}
By an \textbf{isometric action} of a Lie group $G$ on a Riemannian manifold $M$ we mean a smooth and proper homomorphism \mbox{$\Phi:G\to \iso(M)$}. We also denote the action by $\varphi:G\times M\to M,\ (g,q)\mapsto g\cdot q:=\Phi(g)(q)$, or just by $(G,M)$, if no confusion can arise. By this definition, $\varphi$ is smooth and $G\times M\to M\times M, (g,q)\mapsto (g\cdot q, q)$ is proper.
We define the set $M^\reg$ of \textbf{($G$-)regular points} as points on principal orbits. The other points are called \textbf{exceptional} or \textbf{singular}, depending on whether they lie on an exceptional orbit (i.e. it has maximal dimension but its isotropy groups have not the minimal number of connected components) or on a singular one (i.e. it has not maximal dimension).

\begin{definition}\label{d:equiv} We call two $G$-regular points $p,q\in M^\reg$ \textbf{equivalent}, if they can be joined by a broken geodesic whose segments are $G$-transversal geodesics. This means that there is a finite sequence of geodesics $\gamma_0,\dots, \gamma_r:[0,1]\to M^\reg$ sucht that
$$\gamma_0(0)=p,\gamma_r(1)=q,\gamma_i(1)=\gamma_{i+1}(0),\text{ for all } i=0,\dots, r-1$$
and each $\gamma_i$ is perpendicular to one, and hence every, $G$-orbit it meets. It is easily shown that this does indeed define an equivalence relation on $M^\reg$. We denote the equivalence class of $p\in M^\reg$ by $\gothS_p$. We further denote by $\langle\gothS_p\rangle$ the totally geodesic hull of $\gothS_p$, i.e. the connected component containing $p$ of the intersection of all complete totally geodesic submanifolds containing $\gothS_p$.
\end{definition}
Clearly, $\gothS_p\subseteq\langle\gothS_p\rangle\subseteq\Sigma$ holds for every $\Sigma\subseteq M$, which fulfills properties (A)-(C) of a fat section (in \cite{Mag1} such a $\Sigma$ is called \textbf{pre-section}). Furthermore, $\langle\gothS_p\rangle$ intersects with every $G$-orbit. This follows from $\exp_p(\nu_p(G\cdot p))\subseteq\langle\gothS_p\rangle$, and the fact that a minimal geodesic from $p$ to any $G$-orbit already is $G$-transversal. The following construction is crucial for the proof of Theorem \ref{t:reductive_decomp}. It describes how one can obtain a new fat section as a subset out of a given one. This generalizes the construction given in \cite[Section 5.4]{GOT} for orthogonal representations to the case of arbitrary isometric actions of Lie groups on Riemannian manifolds. Our proof modifies the one in loc. cit.

For a vector field $Y$ on $M$ let $\zero(Y)$ denote the set on which $Y$ vanishes. If $Y$ is a Killing field, then the connected components of $\zero(Y)$ are closed and totally geodesic submanifolds of $M$ (see for instance \cite{Kob}). Although we do not assume an isometric action to be effective, we identify the
elements of $\gothg$ with their induced $G$-Killing fields on $M$. We denote the evaluation of the Killing field in $p$ by $X_p$ or $X(p)$ and we put $\phi_g:M\to M,\, q\mapsto g\cdot q$. Let $\Sigma$ be a fat section of $(G,M)$.
For any $X\in\gothg$ the pointwise orthogonal projection \textbf{$\pr_\Sigma X$} of $X|_\Sigma$ to $\Sigma$ is a Killing field of $\Sigma$ (see for instance \cite[Ch. VII, Theorem 8.9]{KN}). Then
$\zero(\pr_\Sigma X)$ consists of those points of $\Sigma$ where $X$ is perpendicular to $\Sigma$. We put
$$\bar\Sigma_p:=(\bigcap_{X\in I_p}\zero(\pr_\Sigma X))^\circ, \text{ where } I_p:=\{X\in\gothg\mid \pr_\Sigma X(p)=0\}.$$
In words, $\bar\Sigma_p$ is the connected component of $p$ of the
common zero set of those projected $G$-Killing fields which vanish
in $p$. It is a closed and totally geodesic submanifold of $\Sigma$.

\begin{proposition}\label{p:induced_section}
Let $\Sigma$ be a fat section of the isometric action $(G,M)$ and let $p\in
\Sigma$ be $G$-regular. Then $\bar\Sigma_p$ is a fat section of $(G,M)$ satisfying
$$\gothS_p\subseteq\,\langle\gothS_p\rangle\subseteq\bar\Sigma_p\subseteq\Sigma.$$
\end{proposition}
\begin{proof}
Let $X\in I_p$ be arbitrary. Then $X$ is a $G$-Killing field which
satisfies $\pr_\Sigma X(p)=0$. By definition, any $q\in\gothS_p$ can
be joined to $p$ by a broken geodesic whose segments are
$G$-transversal geodesics. Since $X_p\,\bot\, T_p\Sigma$, 
we may apply the Lemmas
\ref{l:Killing_Jacobi} and \ref{l:Jacobi_orth} below repeatedly along each
segment to obtain that $X_q\,\bot\,T_q\Sigma$. Thus
$q\in\zero(\pr_\Sigma X)$ and we have shown $\gothS_p\subseteq\zero(\pr_\Sigma X)$. Since $\zero(\pr_\Sigma X)$
is complete and totally geodesic, we also have
$\langle\gothS_p\rangle\subseteq\zero(\pr_\Sigma X)$. Because
$X$ was arbitrary chosen from $I_p$ it follows that
\begin{equation}\label{seq:1}
\langle\gothS_p\rangle\subseteq\bar\Sigma_p.
\end{equation}
Next we show that $\bar\Sigma_p$ is a fat section. Property \ref{pr:1} is obvious
and property \ref{pr:2} follows from
$\langle\gothS_p\rangle\subseteq\bar\Sigma_p$ and the fact that
$\langle\gothS_p\rangle$ intersects every $G$-orbit. Concerning
property \ref{pr:3}, recall that we have to show $\nu_q(G\cdot q)\subseteq T_q\bar\Sigma_p$. We start with the observation:
\begin{equation}
\text{If $q\in\bar\Sigma_p$ is $G$-regular then
$\bar\Sigma_q\subseteq\bar\Sigma_p$. In fact, $I_p\subseteq
I_q$.}
\end{equation}
Using (\ref{seq:1}) and the above observation, we
obtain $\gothS_q\subseteq\bar\Sigma_q\subseteq\bar\Sigma_p$.
Hence
$$\nu_q(G\cdot q)\subseteq
T_q\bar\Sigma_q\subseteq T_q\bar\Sigma_p.$$ In order to show
property \ref{pr:4} we first consider an arbitrary $g\in
N_G(\Sigma)$ and claim:
\begin{eqnarray}
&&g\cdot\zero(\pr_\Sigma X)=\zero(\pr_\Sigma (\Ad_g X)) \text{ and}\\
&&g\cdot\bar\Sigma_p=\bar\Sigma_{g\cdot p}.
\end{eqnarray}
Using that $\Ad_gX(p)=d\phi_g(p)(X_{g^{-1}\cdot p})$, we obtain (3) from the following computation:
\begin{eqnarray*}
g\cdot \zero(\pr_\Sigma X) &=& \{g\cdot p\in\Sigma\mid X_p\,\bot\, T_p\Sigma\} = \{p\in\Sigma\mid X_{g^{-1}\cdot p}\,\bot\, T_{g^{-1}\cdot p}\Sigma\}\\
&=& \{p\in\Sigma\mid d\phi_{g^{-1}}(p)({\Ad_gX}(p))\,\bot\, d\phi_{g^{-1}}(T_p\Sigma)\}\\
&=& \{p\in\Sigma\mid {\Ad_gX}(p)\,\bot\, T_p\Sigma\} = \zero(\pr_\Sigma(\Ad_g X)).
\end{eqnarray*}
From (3) we conclude
$$\Ad_{g^{-1}}X\in I_p \Leftrightarrow X\in I_{g\cdot p}.$$
Now (4) follows from
\begin{eqnarray*}
g\cdot\bar\Sigma_p&=&(\bigcap_{X\in I_p} g\cdot\zero(\pr_\Sigma X))^\circ = (\bigcap_{X\in I_p} \zero(\pr_\Sigma (\Ad_gX)))^\circ\\
&=&(\bigcap_{\Ad_{g^{-1}}X\in I_p} \zero(\pr_\Sigma X))^\circ = (\bigcap_{X\in I_{g\cdot p}} \zero(\pr_\Sigma X))^\circ\\
&=&\bar\Sigma_{g\cdot p}.
\end{eqnarray*}
Let $q\in\bar\Sigma_p$ be an arbitrary $G$-regular point and let $g\in G$ satisfy $g\cdot q\in\bar\Sigma_p$. Then $g\in N_G(\Sigma)$, because $p\in\Sigma$ is $G$-regular. If $q=p$, then (4) and observation (2) imply
$g\cdot\bar\Sigma_p=\bar\Sigma_{g\cdot p}\subseteq\bar\Sigma_p$, and hence $g\cdot\bar\Sigma_p=\bar\Sigma_p$.
For the general case let $\gamma:[0,1]\to M$ be a minimal geodesic from $q$ to $G\cdot p$. Then $\gamma\subseteq\gothS_q$ and there exists some $h\in G$ with $\gamma(1)=h\cdot p\in\bar\Sigma_p$. By the previous arguments $h\cdot\bar\Sigma_p=\bar\Sigma_p$. Note that (1) and (2) show that $\bar\Sigma_x=\bar\Sigma_y$ holds for every $x\in\gothS_y$. Thus $$\bar\Sigma_q=\bar\Sigma_{h\cdot p}=h\cdot\bar\Sigma_p=\bar\Sigma_p.$$
Finally, if $g\cdot q\in\bar\Sigma_p$ then $g\cdot q\in\bar\Sigma_q$ and
$$g\cdot\bar\Sigma_p=g\cdot\bar\Sigma_q=\bar\Sigma_q=\bar\Sigma_p.$$
\end{proof}

\begin{definition}
For a submanifold $N\subseteq M$ an \textbf{$N$-geodesic} $\gamma:[0,\varepsilon)\to M$ is a geodesic of $M$ which emanates perpendicularly from $N$. An \textbf{$N$-Jacobi field} $J$ is a Jacobi field (along an $N$-geodesic $\gamma$), which is induced by a variation of $N$-geodesics.
\end{definition}

Let $\Sigma$ be a fixed fat section through $p\in M^\reg$ and put $N:=G\cdot p$ and $\ke_p:=\nu_p\Sigma$. For $v\in\nu_pN$ we denote by $\gamma_v(t):=\exp_p(tv)$ the $N$-geodesic (which here is the same as a $G$-transversal geodesic) from $p$ in direction of $v$.
The following results, are Lemma 4.3 and Lemma 4.4 of \cite{GOT}. The second lemma characterizes under which conditions an $N$-Jacobi field is perpendicular to a given fat section whereas the first lemma shows that every $N$-Jacobi field, induced by a $G$-Killing field and with the proper initial values, always satisfies the conditions of the second lemma. See also \cite{Mag1}, Section 2.4.
\begin{lemma}[{\cite[Lemma 4.3]{GOT}}]\label{l:Killing_Jacobi}
Let $J$ be an $N$-Jacobi field along $\gamma_v$ with $J(0)\in\ke_p$. If $J$ is the restriction of a $G$-Killing field on $M$ to $\gamma_v$ , then $J$ satisfies $J'(0)+A_v J(0)=0$.
\end{lemma}
\begin{lemma}[{\cite[Lemma 4.4]{GOT}}]\label{l:Jacobi_orth}
Let $J$ be an $N$-Jacobi field along $\gamma_v$ with $J(0)\in\ke_p$. Then $J$ is everywhere perpendicular to $\Sigma$ if and only if $J'(0)+A_v J(0)=0$.
\end{lemma}

From the definition of a fat section it follows that a connected component of the intersection of two fat sections containing a common $G$-regular point is again a fat section. Hence, there is a unique minimal section through every $G$-regular point.
\begin{corollary}\label{c:Killing_perp}
If $\Sigma$ is a minimal section of $(G,M)$, $p\in\Sigma$ is $G$-regular and $X\in\gothg$ is arbitrary, then $X_p\,\bot\, T_p\Sigma$ if and only if $X_q\,\bot\, T_q\Sigma$ holds for all $q\in\Sigma$.
\end{corollary}
\begin{proof}
The above remark and Proposition \ref{p:induced_section} imply that necessarily $\bar\Sigma_p=\Sigma$.
\end{proof}

The corollary is the main ingredient in the
\begin{proof}[Proof of Theorem \ref{t:reductive_decomp}]
Let $p\in\Sigma$ be $G$-regular. By property \ref{pr:3} of a fat section $X_p$ is perpendicular to $\nu_p(G\cdot p)\subseteq T_p\Sigma$. Thus $(\pr_\Sigma X)(p)$ is tangent to $(G\cdot p)\cap\Sigma=N_G(\Sigma)\cdot p$. Let $X_1$ denote an $N_G(\Sigma)$-Killing field which satisfies $X_1(p)=\pr_\Sigma X(p)$. Then $X_1|_\Sigma$ is tangent to $\Sigma$. Putting $X_2:=X-X_1$ we obtain
$$X_2(p)=X(p)-\pr_\Sigma X(p)\in\nu_p(\Sigma).$$
Hence, by Corollary \ref{c:Killing_perp}, it follows that $X_2|_\Sigma$ is always perpendicular to $\Sigma$, and thus $X_1|\Sigma=\pr_\Sigma X$. If $X=Y_1+Y_2$ is another decomposition, then $X_1|_\Sigma=\pr_\Sigma X=Y_1|_\Sigma$. Hence, $0=Y_1|_\Sigma-X_1|_\Sigma$, and it follows that $Z:=Y_1-X_1$ is a $Z_G(\Sigma)$-Killing field, because the latter are characterized as those $G$-Killing fields which vanish everywhere on $\Sigma$. Since
$$X=X_1+X_2=Y_1+Y_2=X_1+Z+Y_2,$$
it follows that $Y_2=X_2-Z$, and clearly $Z$ is uniquely determined.

We next show the inclusion $\Ad_g(\gothm)\in\gothm$. Let $X\in\gothm,\, g\in N_G(\Sigma)$ and $p\in\Sigma^\reg$ be arbitrary elements (we note here that due to \cite[Lemma 2.3.1]{Mag1} the sets of $G$-regular, resp. $W$-regular points on $\Sigma$ coincide). Then $g\cdot p\in\Sigma^\reg$ and it follows that
$$(\Ad_g X)(g\cdot p)=d\phi_g(p)(X_p)\,\bot\, d\phi_g(p)(T_p\Sigma)=T_{g\cdot p}\Sigma.$$
Thus $(\Ad_g X)(g\cdot p)\in\nu_{g\cdot p}\Sigma$ and Corollary \ref{c:Killing_perp} yields $\Ad_g X\in\gothm$. The inclusion $[\gothn,\gothm]\subseteq\gothm$ of course follows from $\Ad_g(\gothm)\subseteq\gothm$.
\end{proof}

\begin{corollary}\label{c:G-W-invariant1} Let $\Sigma, N, H$ and $W$ as in Theorem \ref{t:reductive_decomp}. Then
\begin{enumerate}
\item $G\!/\!H$ admits an adapted Riemannian metric if and only if $\gothm/\gothh$ carries an $\Ad_G(N)$-invariant scalar product and $W$ is covered by the product of a compact Lie group and a vector group.
\item If $N$ is compact, then $G\!/\!H$ admits an adapted Riemannian metric.
\end{enumerate}
\end{corollary}
\begin{proof}
Clearly, the first statement implies the second one. By Theorem \ref{t:reductive_decomp}, the decomposition $\gothg/\gothh=\gothn/\gothh\oplus\gothm/\gothh$ is direct and $\Ad_G(N)$-invariant. If $\langle\cdot,\cdot\rangle$ on $\gothg/\gothh$ is induced by an adapted $(G$-$W)$-invariant Riemannian metric on $G\!/\!H$, then its restriction to $\gothm/\gothh$, resp. $\gothn/\gothh$ yield an $\Ad_G(N)$-invariant inner product on $\gothm/\gothh$, resp. an $\Ad_W$-invariant scalar product on $\gothn/\gothh$. Since $W$ is a Lie group, the existence of such a scalar product is equivalent to the existence of a bi-invariant Riemannian metric on $W$. Using \cite[Proposition 3.34]{CE} it follows that $W$ is covered in the way stated in the corollary. Conversely, we may patch a pair of invariant scalar products on $\gothm/\gothh$ and $\gothn/\gothh$ together to form an $\Ad_G(N)$-invariant inner product on $\gothg/\gothh$ such that both factors are perpendicular to each other. This may be extended to a left-$G$-invariant metric on $G\!/\!H$, which is easily seen to be right-$W$-invariant.
\end{proof}

Note that in the case that the action is either polar (i.e. $\gothh=\gothn$) or has trivial copolarity (i.e. $\gothh=\gothm$) any left-invariant metric on $G\!/\!H$ is adapted to $\Sigma$.
With respect to a $(G$-$W)$-invariant Riemannian metric $W$ is totally geodesic in $G\!/\!H$. If the metric is adapted and if $G\!/\!N$ carries the metric induced from $G\!/\!H\twoheadrightarrow G\!/\!N$, then $\gothm/\gothh$ is canonically isometric to $\gothg/\gothn$.

\section{Weyl's Integration formula}
The main tool in this section is a generalization of Fubini's Theorem for general surjective submersions. The case of a Riemannian submersion is already well known in the literature (cf. \cite[Chapter II, Theorem 5.6]{Sa}). In this case integration of a function on the total space with respect to the Riemannian measure can be reduced to an integration over the fibre and the base space. In our case, integration over the fibre will be weighted by a factor, which equals one for a Riemannian submersion. Since the proof is just a slight modification of the one in loc. cit., we omit it here.

In the following we denote the Riemannian measure on a Riemannian manifold $(M,g)$ by $\mu_g$. Let $\pi:(M,g)\to(N,h)$ be a surjective submersion between two Riemannian manifolds. For $q\in N$ let $\delta_q$ denote the real function
$$\delta_q:\pi^{-1}(q)\to\R,\, p\mapsto |\det(d\pi(p)|_{\kh_p})^{-1}|,$$
where $\kh_p$ denotes the horizontal space to the fibre
$\pi^{-1}(q)$ in $T_pM$. Recall that $|\det f|$ for a linear
homomorphism $f:X\to Y$ between two equal dimensional Euclidean
spaces is defined via the usual determinant as $|\det f|:=|\det
A\circ f|$, where $A:Y\to X$ is an arbitrary auxiliary linear
isometry. If $X$ and $Y$ have unequal dimension, we set $|\det
f|:=0$. Clearly, this is well defined as it is independent of the
choice of $A$. For Riemannian submersions
of course $\delta_q\equiv 1$, because each horizontal space
$\kh_p$ along the fibre over $q$ is mapped isometrically onto
$T_qN$. In general this need not be the case, although we always
have $\delta_p>0$. Using local frames it is not difficult to show
that $\delta_p$ is smooth along $\pi^{-1}(q)$. For a function
$f:M\to\R$ we denote its restriction to the fibre by
$f_q:=f|_{\pi^{-1}(q)}$ and we denote the Riemannian metric on
$\pi^{-1}(q)$ induced by $g$ with $g_q$. If $f_q$ is integrable with
respect to the weighted measure $\delta_q\mu_{g_q}$ on $\pi^{-1}(q)$
we put
$$\bar f(q):=\int_{\pi^{-1}(q)} f_q\delta_q\mu_{g_q}.$$
\begin{proposition}[Fubini's Theorem for submersions]\label{p:Fubini}
Let $\pi:(M,g)\twoheadrightarrow (N,h)$ be a surjective submersion between two Riemannian manifolds. If $f\in\kc_c(M)$ (resp. $f$ is integrable on $M$), then $\bar f\in\kc_c(N)$ (resp. $f_q$ is integrable for almost all $q\in N$ and $\bar f$ is integrable on $N$). Furthermore,
$$\int_M f\, d\mu_g=\int_N\bar f\,d\mu_h=\int_N\left(\int_{\pi^{-1}(q)}f_q\, \delta_q d\mu_{g_q}\right)d\mu_h.$$
\end{proposition}

One feature of a minimal section $\Sigma$ of the isometric action $(G,M)$ is that the principal isotropy groups are constant along $\Sigma$. I.e. $G_p=H$ for all $p\in\Sigma^\reg$. We now describe how for each $s\in\Sigma$ the volume of the orbit $G\cdot s$ can be computed by the orbit map $\omega_s:G\!/\!H\to G\cdot s$. As mentioned in the introduction, this map is a diffeomorphism, a covering or a fibre bundle, depending on whether $s$ is regular, exceptional or singular. The covering in the second case is $|G_s\!/\!H|$-fold. We first assume that $s$ is regular. Applying the transformation formula yields
$$\vol(G\cdot s)=\int_{G\cdot s} 1\, dy=\int_{G\!/\!H}|\det(d\omega_s(gH))|\, d(gH).$$
By assumption, the metric on $G\!/\!H$ is left-$G$-invariant. Since
$\omega_s$ is equivariant with respect to the $G$-actions on $G\!/\!H$
and $G\cdot s$ and isometries have a determinant with absolute value
$1$, it follows that
$$|\det d\omega_s(gH)|=|\det d\omega_s(eH)|.$$
This yields the formula
$$\vol(G\cdot s)=\vol(G\!/\!H) |\det d\omega_s(eH)|.$$
This shows that $\vol(G\!/\!H)$ is finite if $\vol(G\cdot s)$ is finite. By a similar argument one shows that the last statement is actually an ``if and only if''-statement. If $s$ is exceptional, then $G_s\!/\!H$ is a finite group since $G_s$ is compact by the properness of the action and $H$ is open in $G_s$. Fubini's Theorem applied to the $|G_s\!/\!H|$-fold covering $G\!/\!H\to G\!/\!G_s$ gives
\begin{lemma}
$\vol(G\!/\!H)$ is finite, if and only if $\vol(G\cdot s)$ is finite. Furthermore
$$\vol(G\cdot s)=\tfrac{1}{|G_s\!/\!H|}\, \vol(G\!/\!H)\, |\det d\omega_s(eH)|.$$
\end{lemma}
Applying Fubini's Theorem to the Riemannian submersion
$W\hookrightarrow G\!/\!H\twoheadrightarrow G\!/\!N$ yields that $W$ and
$G\!/\!N$ have finite volume if $G\!/\!H$ has. As the action of $W$ is free on $\Sigma^\reg$ the above lemma implies
$$\vol(W\cdot s)=\tfrac{1}{|W_s|}\, \vol(W)\, \delta_\kd(s)$$
for all non-singular $s\in\Sigma$. Here we have denoted the determinant expression by $\delta_\kd$ in order to avoid confusion. A closer look on the map
$$d\omega_s(eH):\gothg/\gothh\to T_s(G\cdot s),\ X+\gothh\mapsto X_s$$
reveals that it maps the parts $\gothn/\gothh$ and $\gothm/\gothh$, of the decomposition $\gothg/\gothh=\gothn/\gothh\oplus\gothm/\gothh$ in Theorem \ref{t:reductive_decomp}, onto $T_s(W\cdot s)$, resp. onto the normal space $\nu_s\Sigma$ of $T_s\Sigma$ in $T_s M$. Thus, $d\omega_s(eH)$ has block diagonal form. Recalling the definition of $\delta_\ke$ in the introduction:
$$|\det d\omega_s(eH)|=\delta_\kd(s)\, \delta_\ke(s).$$
From these considerations, the proof of Theorem \ref{t:integration} (iv) follows.

The proof of formula (i) will also be an application of Fubini's
Theorem. For this purpose consider the map $\pi:M^\reg\to G\!/\!N$,
which is defined via property (D) of a fat section. Namely, for any
$p\in M^\reg$, there is a unique translate $g\cdot\Sigma$ of
$\Sigma$ containing $p$. The coset $gN$, where
$N=N_G(\Sigma)$ is the normalizer of $\Sigma$ in $G$, is uniquely
determined by $p$ and we put $\pi(p):=gN$. Alternatively, consider
the bundle $G\times_N\Sigma^\reg$ associated with the $N$-principal
bundle $G\twoheadrightarrow G\!/\!N$. The action map $\varphi:G\times
M\to M$ induces a diffeomorphism
$\tilde\varphi:G\times_N\Sigma^\reg\to M^\reg$. Now the natural
projection map $G\times_N\Sigma^\reg\twoheadrightarrow G\!/\!N$ has
fibres which under $\tilde\varphi$ are the $G$-translates of
$\Sigma$. In this way we see immediately that $\pi$ is a smooth
surjective submersion. The horizontal space to the fibre
$g\cdot\Sigma$ is given by
$$\nu_{g\cdot s}(g\cdot\Sigma)=d\phi_g(s)(\nu_s\Sigma).$$
Since $\pi$ is equivariant with respect to the isometric actions of
$G$ on $M^\reg$ and $G\!/\!N$, it follows that
$\delta:g\cdot\Sigma^\reg\to\R,\, p\mapsto |\det(d\pi(p)|_{\nu_p
(g\cdot\Sigma)})^{-1}|$ is $G$-invariant and hence completely
determined by its values on $\Sigma^\reg$. Its relation to
$\delta_\ke$ is the following:
\begin{lemma}
$\delta_\ke(s)=\delta(s)$ for all $s\in\Sigma^\reg$.
\end{lemma}
\begin{proof}
Since the metric on $G\!/\!H$ is adapted, the spaces $\gothg/\gothn$ and $\gothm/\gothh$ are isometric by
\begin{equation}\label{e:isometric}
\gothm/\gothh\to\gothg/\gothn,\, X+\gothh\mapsto X+\gothn.
\end{equation}
This in turn is induced by the differential of the $W$-principal bundle $G\!/\!H\twoheadrightarrow G\!/\!N$. With this identification we claim
\begin{equation}\label{e:0}
d\omega_s(eH)|_{\gothm/\gothh}=(d\pi(s)|_{\nu_s\Sigma})^{-1}.
\end{equation}
In order to compute $d\pi(s)$ we first note that property (C) of a fat section implies $\nu_s\Sigma\subseteq T_s(G\cdot s)$. Hence, for every $v\in\nu_s\Sigma$ there is a one parameter group of $G$ induced by some $X\in\gothg$ such that $\frac{d}{dt}|_{t=0}\exp(tX)\cdot s=v$. Accordingly
$$d\pi(s)(v)=\left.\frac{d}{dt}\right|_{t=0}\pi(\exp(tX)\cdot s)=\left.\frac{d}{dt}\right|_{t=0}\exp(tX)N=X+\gothn.$$
By the isometry in (\ref{e:isometric}) we may assume $X\in\gothm$. Hence
$$d\omega_s(eH)\circ d\pi(s)(v)=d\omega_s(eH)(X+\gothm)=X_s=v.$$
From this the claimed identity (\ref{e:0}) follows and thus $\delta_\ke=\delta$.
\end{proof}
In the following we tacitly use the fact that $\Sigma^\reg$ and $M^\reg$ lie dense in $\Sigma$, resp. $M$. If we apply Fubini's Theorem to $\pi$, we get for all $f\in L^1(M)$:
\begin{equation}\label{e:1}
\int_M f(x) dx=\int_{G\!/\!N}\bar f(gN)d(gN),
\end{equation}
where
\begin{equation}\label{e:2}
\bar f(gN)=\int_{g\cdot\Sigma}f(x)\delta(x)dx.
\end{equation}
Applying the transformation formula to (\ref{e:2}) with respect to the isometry $\phi_g$ we obtain:
\begin{equation}\label{e:3}
\bar f(gN)=\int_{\Sigma}f(g\cdot s)\delta(g\cdot s)ds=\int_{\Sigma}f(g\cdot s)\delta(s)=\int_{\Sigma}f(g\cdot s)\delta_\ke(s)ds.
\end{equation}
Inserting (\ref{e:3}) into (\ref{e:1}), yields formula (i) of Theorem \ref{t:integration}. This in turn implies formula (ii).

For the proof of (iii) we need the invariance property of $\delta_\ke$ from Proposition \ref{p:delta}.
\begin{proof}[Proof of Proposition \ref{p:delta}]
Fix an orthonormal frame $X_1+\gothh,\dots,X_m+\gothh$ of
$\gothm/\gothh$ and a local orthonormal frame $Y_1,\dots,
Y_m\in\Gamma(\Sigma,\nu(\Sigma))$ of the normal bundle $\nu(\Sigma)$
of $\Sigma$ on a neighborhood of $s\in\Sigma$. We define the
auxiliar pointwise linear isometry
$$A_s:\nu_s\Sigma\to\gothm/\gothh$$
by mapping $Y_i(s)$ to $X_i+\gothh$. Then $A_s$ varies smoothly in $s$. Also note, that $d\omega_s(eH)$ depends smoothly on $s$, because the action $\varphi$ is smooth. The usual determinant of square matrices is a polynomial in the matrix entries and thus
\begin{equation}\label{e:4}
\det(A_s\circ d\omega_s(eH)|_{\gothm/\gothh})
\end{equation}
is smooth in $s$ as a composition of smooth maps. For
$s\in\Sigma^\reg$ expression (\ref{e:4}) does not vanish on a small
neighborhood of $s$. Hence its sign does not change there and thus
$$\delta_\ke(s)=|\det(A_s\circ d\omega_s(eH)|_{\gothm/\gothh})|$$
is smooth on $\Sigma^\reg$. If $s$ is not $G$-regular, expression (\ref{e:4}) may vanish and the sign may change on a neighborhood of $s$. Hence $\delta_\ke$ need not be smooth on the whole of $\Sigma$, but it is still continuous.

Since $H\subseteq N$ is a normal subgroup, $\omega_s$ is also
equivariant with respect to the right action of $W=N\!/\!H$ on $G\!/\!H$ and
the left action of $W$ on $W\cdot s$. Since the metric on $G\!/\!H$ is
also right $W$-invariant, the $W$-action on $G\!/\!H$ is isometric. From
this the $W$-invariance of $\delta_\ke$ follows. As a general
fact, any $W$-invariant continuous function on $\Sigma$ has a
natural $G$-invariant continuous extension to $M$. If it is smooth on $\Sigma^\reg$, the extension is smooth on
$M^\reg$ (cf. \cite[Corollary 2.1.2 and Proposition 2.7.1]{Mag1}).
\end{proof}

\begin{proof}[Continuation of the proof of Theorem \ref{t:integration} (iii)] Formula (ii) and Proposition \ref{p:delta} imply
$$\sqrt[p]{\vol(G\!/\!N)\delta_\ke}\, f|_\Sigma\in L^p(\Sigma)^W.$$
Thus $\Theta_p$ is well defined. It is obviously a linear map.
Formula (ii) also implies that $\Theta_p$ is an isometry and hence
injective and continuous. It remains to prove the surjectivity of
$\Theta_p$. Let $\tilde f\in L^p(\Sigma)^W$ be arbitrary. Since
$\Sigma^\reg$ is dense in $\Sigma$, we can approximate $\tilde f$ by
a sequence $\tilde f_n\in\kc_c(\Sigma^\reg)^W$ with respect to the
$L^p$-norm. We consider its unique continuation
$h_n\in\kc(M^\reg)^G$. We claim
that $\supp(h_n)=G\cdot\supp(\tilde f_n)$ is compact. Let therefore
$x_k=g_k\cdot y_k\in\supp(h_n)$ with $g_k\in G$ and
$y_k\in\supp(h_n)$ be an arbitrary sequence for a fixed
$n\in\N.$ 
Since $G\!/\!N$ is assumed to be compact, there is some $g\in G$ such
that, after passing to a subsequence if necessary,
$\lim_{k\to\infty}g_kN=gN$. Hence there is a sequence $m_k\in N$
with $\lim_{k\to\infty}g_km_k=g$. Since $\supp(\tilde f_n)$ is
compact and invariant under $N$, we may assume again that
$\lim_{k\to\infty}m_k^{-1}\cdot y_k=y\in\supp(\tilde f_n)$. It thus
follows that $x_n=(g_km_k)\cdot(m_k^{-1}\cdot y_k)$ has a convergent
subsequence. This shows that $\supp(h_n)$ is compact for every $n\in\N$. Let now
$$f_n:=\frac{h_n}{\sqrt[p]{\vol(G\!/\!N)\delta_{\ke}}}\in\kc_c(M^\reg)G.$$
By the above considerations this is well defined and satisfies
$\Theta_p(f_n)=\tilde f_n$. Using formula (ii) again we see that
$f_n$ is a Cauchy sequence and thus converges to some $f\in
L^p(M)^G$. By continuity of $\Theta_p$ we conclude
$\Theta_p(f)=\tilde f$, which proves the surjectivity of $\Theta_p$.
\end{proof}

\section{The Direct Sums of the Standard Representations $\rho_n, \mu_n$ and $\nu_n$}
\begin{proof}[Proof of Proposition \ref{p:direct_sum}]
The first step of the proof that the given $\Sigma$ is indeed a minimal section is to construct the equivalence class $\gothS_v$ (see Definition \ref{d:equiv}) through a suitable $G$-regular point $v\in\Sigma$. We then show that $\Sigma$ is already the $\R$-linear span of $\gothS_v$. The copolarity is then of course the difference between the dimension of $\Sigma$ and the cohomogeneity of the representation.
For $2\le k\le n-1$ let
$$v=(e_1,\dots,e_k)=:\begin{pmatrix}\textbf{1}\\ \textbf{0}\end{pmatrix}.$$
(For $k\ge n$ the proof is basically the same). Then $v\in\Sigma$ and it is a $G$-regular point with isotropy group equal to $H$ in the table above. From this we can compute the cohomogeneity of the representation $\varphi$. Let $\ka(\K)$ denote the set of all skew-symmetric, skew-hermitian, resp. skew-quaternionic-hermitian $k\!\times\! k$-matrices and let $\ks(\K)$ denote the set of all symmetric, hermitian, resp. quaternionic-hermitian $k\!\times\! k$-matrices. Then the tangent space at the orbit through $v$ is described by:
$$\gothg\cdot v=\left\{\begin{pmatrix}A\\B\end{pmatrix}\mid A\in\ka(\K),\ B\in\K^{(n-k)\!\times\! k} \right\}.$$
The (real) inner product of two elements $x,y\in V$ is given by
$$\langle x|y\rangle=\Re(\tr(x^*y)),$$
where $*$ denotes transposition, complex-conjugate transposition or quaternionic-con\-ju\-ga\-te transposition. 
The normal space of the $G$-orbit \mbox{through $v$ is:}
$$\nu_v(G\cdot v)=\left\{\begin{pmatrix}C\\\textbf{0}\end{pmatrix}\mid C\in\ks(\K)\right\}.$$

The point $w=(e_1,2e_2,\dots,ke_k)$ is $G$-regular and contained in $\nu_v(G\cdot v)$. For $\K=\R$
$$\gothg\cdot v\cap\gothg\cdot w=\left\{\begin{pmatrix}\textbf{0}\\B\end{pmatrix}\mid B\in\R^{(n-k)\!\times\! k}\right\}.$$
This coincides with the normal space of $\Sigma$ in $V$. Therefore, in this case the formula
$$(U\cap W)^\bot=U^\bot+W^\bot$$
for subspaces $U,W\subseteq V$ already implies
$$\Sigma=\nu_v(G\cdot v)+\nu_w(G\cdot w)$$
and we are done with the proof. However, in the case of $\K=\C$ or $\K=\H$ we have
$$\gothg\cdot v\cap\gothg\cdot w=\left\{\begin{pmatrix}D\\B\end{pmatrix}\mid D\in\Pu_k(\K),\ B\in\R^{(n-k)\!\times\! k}\right\},$$
where $\Pu_k(\K)$ is the set of $(k\!\times\! k)$-diagonal matrices with entries in the imaginary numbers, resp. pure quaternions.
In both cases let $u\in\nu_v(G\cdot v)$ be the $(n\!\times\! k)$-matrix
$$u_{ij}:=\left\{\begin{array}{cl}0 & \text{if } i=j \text{ or } i>k,\\
1 & \text{else.}\end{array}\right.$$
It is clear that $u$ is $G$-regular, since its rank is equal to $k$.
Now
\begin{equation}\tag{$*$}
\gothg\cdot v\cap\gothg\cdot w\cap\gothg\cdot u=\left\{\begin{pmatrix}\begin{array}{ccc}\lambda&&0\\ &\ddots&\\ 0&&\lambda\\\end{array}\\B\end{pmatrix}\mid \lambda\in\Pu(\K),\ B\in\R^{(n-k)\!\times\! k}\right\}.\end{equation}
In fact, let $a\in\gothg$ be such that $a\cdot u\in \gothg\cdot v\cap\gothg\cdot w$. Then for all $i,j$:
\begin{equation}\tag{$**$}
(a\cdot u)_{ij}=\sum_{l=1}^n
a_{il}u_{lj}=\sum_{\tiny\begin{array}{c}l=1\\l\ne
j\end{array}}^k a_{il}.\end{equation} For $i\ne j$ in
$\{1,\dots,k\}$ this expression vanishes and thus for $i\in\{1,\dots, k\}$ fixed:
$$0=\sum_{\tiny\begin{array}{c}l=1\\l\ne j\end{array}}^k a_{il}-\sum_{\tiny\begin{array}{c}l=1\\l\ne m\end{array}}^k a_{il}=a_{im}-a_{ij}$$
holds for all $j,m\in\{1,\dots, k\}-\{i\}$. Thus
$(k-2)a_{ij}=-a_{ii}\in\Pu(\K)$, which implies $a_{ij}=-\bar a_{ij}$
for all $i\ne j$ (in case of $k=2$ this follows from $(**)$ and
$a\cdot u\in\gothg\cdot v\cap\gothg\cdot w$). Together with
$a_{ij}=-\bar a_{ji}$ we conclude that $a_{ij}=a_{ji}$ holds for all
$i\ne j$, and finally
$$(a\cdot u)_{ii}=\sum_{\tiny\begin{array}{c}l=1\\l\ne i\end{array}}^k a_{il}=(n-1)a_{12}$$
holds for all $i=1,\dots, k$, which proves equation $(*)$.

The last element we consider is $\tilde u:=(2e_2,e_1,3e_3,\dots,ke_k)\in\nu_w(G\cdot w)$. This is again a $G$-regular point in $\gothS_v$, and it is easily verified that
$$\gothg\cdot v\cap\gothg\cdot w\cap\gothg\cdot u\cap\gothg\cdot\tilde u=\left\{\begin{pmatrix}\textbf{0}\\B\end{pmatrix}\mid B\in\R^{(n-k)\!\times\! k}\right\}.$$
As in the case of $\K=\R$ it now follows that
$$\Sigma=\nu_v(G\cdot v)+\nu_w(G\cdot w)+\nu_u(G\cdot u)+\nu_{\tilde u}(G\cdot \tilde u).$$

Concerning the normalizer, centralizer and the fat Weyl group of $\Sigma$ in $G$ we have:
\begin{eqnarray*}
N(\Sigma)&=&S(\On(k)\!\times\!\On(n-k)),\ S(\Un(k)\!\times\!\Un(n-k)), \text{ resp. } \Spn(k)\!\times\!\Spn(n-k),\\
Z(\Sigma)&=&\{\textbf{1}\}\!\times\!\Son(n-k),\ \{\textbf{1}\}\!\times\!\Sun(n-k), \text{ resp. } \{\textbf{1}\}\!\times\!\Spn(n-k),\\
W(\Sigma)&=&\On(k),\ \Un(k), \text{ resp. }\Spn(k).
\end{eqnarray*}
\end{proof}
As a remark, due to the following lemma each minimal section appearing in Proposition \ref{p:direct_sum} can be written as the fixed point set of $\tilde H$, where $\tilde H$ is a principal isotropy group for a representation of a larger group $\tilde G$, which has the same orbits as $G$.
\begin{lemma} Let $n\ge 2$ be arbitrary and let $1\le k\le n-1$. Then
\begin{enumerate}
\item The representation $k\cdot\rho_n$ of $\Son(n)$ on $k$ copies of
$\R^n$ has the same orbits as the corresponding representation of
$\On(n)$ on $k$ copies of $\R^n$.
\item The representation $k\cdot\mu_n$ of $\Sun(n)$ on $k$ copies of
$\C^n$ has the same orbits as the corresponding representation
of $\Un(n)$ on $k$ copies of $\C^n$.
\end{enumerate}
\end{lemma}
In Straume \cite{S} (see also \cite{GOT}) it is asked if every minimal section can always be described as a component of the fixed point set of a principal isotropy group, if one passes to the largest (possibly non-connected) Lie group inducing the same orbits.

\begin{proof}[Continuation of the proof of Proposition \ref{p:direct_sum}]
We still have to compute $\delta_\ke$. The scalar product $\langle x|y\rangle=\Re(\tr(x^*y))$ on $V$ also defines a scalar product on $\gothg$. It induces an adapted $(G$-$W)$-invariant Riemannian metric on $G\!/\!H$. The elements of $\gothn$, resp. $\gothn^\bot$ are block matrices of the form
$$\begin{pmatrix}A&\mathbf{0}\\ \mathbf{0}&D\end{pmatrix},\text{ resp. } \begin{pmatrix}\mathbf{0}&E^*\\ E & \mathbf{0}\end{pmatrix},$$
where $A\in\on(k), \un(k)$, resp. $\spn(k)$, $D\in\on(n-k), \un(n-k)$, resp. $\spn(n-k)$ and $\tr(A)+\tr(D)=0$, and $E\in\K^{(n-k)k}$ is an arbitrary matrix. $\gothn^\bot$ is isometric to $\gothm/\gothh$ and the orbit map $d\omega_p(e):\gothn^\bot\to\nu_p\Sigma$ is just matrix multiplication of $p$ by elements of $\gothn^\bot$ from the left
$$d\omega_p(e)\begin{pmatrix}\mathbf{0}&E^*\\ E & \mathbf{0}\end{pmatrix}=\begin{pmatrix}\mathbf{0}&E^*\\ E & \mathbf{0}\end{pmatrix}\begin{pmatrix}B^t\\\textbf{0}\end{pmatrix}=\begin{pmatrix}\textbf{0}\\ EB^t\end{pmatrix}.$$
For $i=k+1,\dots,n$ and $j=1,\dots,k$ let $E_{ij}\in\R^{(n-k)k}$ denote the corresponding elementary matrix. Then
$$e_{ij}:=\begin{pmatrix}\textbf{0} & \frac{-1}{\sqrt2}(E_{ij})^*\\\frac{1}{\sqrt2}E_{ij} & \textbf{0}\end{pmatrix} \text{ and } f_{ij}:=\begin{pmatrix}\textbf{0}\\ E_{ij}\end{pmatrix}$$
form ON-bases of $\gothn^\bot$, resp. $\nu_p\Sigma$, if we view both as $\K$-vector spaces. Identifying $f_{ij}$ with $e_{ij}$ yields a linear isometry between $\nu_p\Sigma$ and $\gothn^\bot$. We arrange the $e_{ij}$ as follows:
$$e_{k+1,1},\dots, e_{k+1,k},e_{k+2,1},\dots, e_{k+2,k},\dots,e_{n,1},\dots,e_{n,k}.$$
With respect to this ordering, $d\omega_p(e)$ has the matrix representation $\diag(\frac{1}{\sqrt2}B,\dots,\frac{1}{\sqrt2}B)$,
the block $\frac{1}{\sqrt2}B$ appearing $(n-k)$-times. We are interested in the absolute value of the determinant of $d\omega_p(e)$ as a map between real vector spaces. In \cite{A} one can find the following relation for the determinant $\det(P)$ of a complex, resp. quaternionic matrix $P$ and the determinant of its realification $\det{_\R}(P)$:
$$\det{_\R}(P)=|\det(P)|^d.$$
Since determinants are multiplicative, even if $K=\H$, the claimed formula follows.
\end{proof}

\section*{Acknowledgement}The Author received financial support by the DFG-Schwerpunkt 1145 \glqq Globale Dif\-fer\-en\-tial\-geo\-metrie\grqq\ and the University of Münster. This article is based on parts of the authors doctoral thesis \cite{Mag1}.
\bibliography{Literaturverzeichnis}
\bibliographystyle{alpha}
\end{document}